\documentclass [12pt]{amsart}

\usepackage{latexsym,amssymb, amsmath}

\usepackage{simplemargins}
\setleftmargin{2.75cm}
\setrightmargin{2.75cm}
\settopmargin{2.0in}
\setbottommargin{2.0in}
\setallmargins{2.75cm}

\def\pf{\noindent\emph{Proof: }}
\def\stop{\hfill$\Box$}

\usepackage{amsmath, amsfonts, amssymb, amsthm}
\newtheorem{thm}{Theorem}

\newtheorem{lemma}[thm]{Lemma}

\numberwithin{thm}{section}
\DeclareMathOperator{\Tr}{Tr}

\begin{document}

\title 
[Uniform Estimates for the Constant Mean Curvature Solutions ] 
{Uniform Estimates for the Constant Mean Curvature Solutions of the Vacuum Einstein Constraint Equations on Compact Manifolds  }

\author{Yu Yan }
\address {Department of Mathematics\\The University of British Columbia  \\Vancouver, B.C., V6T 1Z2\\ Canada}

\email{yyan@math.ubc.ca}


\begin{abstract}
We give some uniform estimates for constant mean curvature solutions of the conformal vacuum Einstein constraint equations on compact manifolds.  Existence of those solutions was given in \cite{isen1}.

\end{abstract}

\maketitle

\section {Introduction}

Let $M^3$ be a 3-dimensional Riemannian manifold with metric $\gamma$, let $K$ be a symmetric $(0,2)$ tensor on $M$.  It was shown by Y. Choquet-Bruhat \cite{FBru} that any initial data set $(M^3, \gamma, K)$ can be viewed as an embedded hypersurface in a spacetime with $\gamma$ as the induced metric and $K$ as its second fundamental form, provided it satisfies the following {\it Einstein constraint equations}
\begin{eqnarray}
\label{eq:einstein}
R^M+(\Tr _{\gamma}K)^2-\|K\|^2 & = & 16 \pi T_{nn},\\
K^j_{i;j}-K^j_{j;i} & = & 8 \pi T_{ni}, 
\end{eqnarray}
where $R^M$ is the scalar curvature of $M$ and $T_{ab}$ is the stress-energy tensor that describes the matter content of the ambient spacetime.  Therefore, constructing an initial data set that satisfies those constraint equations is the first step in understanding the global evolution of the spacetime.

The most widely used approach to the constraint equations in the vacuum case ($T_{ab}=0$) is the {\it conformal method} \cite {BarIsen}, which divides the initial data on $M^3$ into the ``Free (conformal) Data'' (a Riemannian metric $\lambda _{ij}$, a divergence-free, trace-free symmetric tenser $\sigma _{ij}$, and a scalar function $\tau$), and the ``Determined Data''(a strictly positive scalar function $\phi$ and a vector field $W^i$). Given a choice of the conformal data, if one can find a determined data satisfying the determined elliptic PDE system
\begin{eqnarray}
\label{eq:LW}
\label{eq:phi}
\nabla _i(LW)^i_j & = & \frac{2}{3} \phi ^6 \nabla _j \tau ,\\
\label{eq:phi}
\Delta \phi & = & \frac{1}{8} R \phi - \frac{1}{8}(\sigma^{ij}+LW^{ij})(\sigma_{ij}+LW_{ij})\phi^{-7}+\frac{1}{12}\tau^2\phi^5,
\end{eqnarray}
then the two sets of data can be combined to produce an initial data set which satisfies the constraint equations by
\begin{eqnarray}
\label{eq:gamma}
\gamma_{ij} & = & \phi ^4 \lambda _{ij},\\
\label{eq:K}
K_{ij} & = & \phi ^{-2}(\sigma_{ij}+LW_{ij})+\frac{1}{3}\phi ^4 \lambda _{ij}\tau.
\end{eqnarray}

This conformal method has been very successful in finding constant mean curvature solutions ($\tau = $ constant) of the vacuum constraint equations. With $\tau$ constant, (\ref{eq:LW}) implies $LW=0$, and consequently (\ref{eq:phi}) becomes the Lichnerowicz equation
\begin{equation}
\label{eq:lich}
\Delta \phi  =  \frac{1}{8} R \phi - \frac{1}{8}\|\sigma\|^2\phi^{-7}+\frac{1}{12}\tau^2\phi^5.
\end{equation}
This equation is conformally invariant in the sense that it has a solution with respect to conformal data ($\lambda, \sigma, \tau$) if and only if it has a solution with respect to conformal data ($\psi^4\lambda, \psi^{-2}\sigma, \tau$) for some function $\psi >0$.  

The Yamabe invariant of $\lambda$ is defined as 
\begin{eqnarray*}
Y(M, \lambda)& = &\inf \left\{ \frac{\int_M R_{\bar {\lambda}} dv_{\bar {\lambda}}
}{(\int_M dv_{\bar {\lambda}} )^\frac{n-2}{n} }: \text{  } \bar{\lambda}=\psi(x)^{\frac{4}{n-2}}\lambda, \,\, \psi (x)>0, \,\,\psi \in W^{1,2}(M) \right\} \\
&  =  & \inf \left\{ \frac{\int_M (|\nabla \psi|^2 + \frac{n-2}{4(n-1)}
R(\lambda)\psi ^2)dv_{\lambda} }{(\int_M \psi^{\frac{2n}{n-2}}dv_{\lambda})^{\frac{n-2}{n}} }: \text{  } \psi(x)>0, \psi \in W^{1,2}(M) \right\}.
\end{eqnarray*}

\noindent
By the Yamabe Theorem (\cite{Au}, \cite{S1}, \cite{Tr}), $Y(M,\lambda)$ being positive, zero, and negative implies that $\lambda$ is conformal to a metric of positive constant, zero, and negative constant scalar curvature, respectively.  Therefore, all Riemannian metrics on $M^3$ can be divided into three classes by their Yamabe invariant: Yamabe positive $\mathcal{Y}^+(M)$, Yamabe zero $\mathcal{Y}^0(M)$, and Yamabe negative $\mathcal{Y}^-(M)$. Since we can also divide $\sigma$ and $\tau$ into $\|\sigma\|^2 \not\equiv 0$, $\|\sigma\|^2 \equiv 0$ and $\tau \neq 0$, $\tau =0$, there are totally twelve classes of conformal data to consider.  In \cite{isen1} Isenberg gave a theorem which completely determines for which of the twelve classes equation (\ref {eq:lich}) can be solved by the conformal method.  His theorem can be organized into the following table

\vspace{.05in}
\begin{center}
\begin{tabular}{c|c|c|c|c}
 &  $\|\sigma\|^2 \equiv 0, \tau =0$  & $\|\sigma\|^2 \equiv 0, \tau \neq0$  & $\|\sigma\|^2 \not\equiv 0, \tau =0$  & $\|\sigma\|^2 \not\equiv 0, \tau \neq0$  \\
\hline
$\mathcal{Y}^+(M)$ & No & No & Yes & Yes\\
\hline
$\mathcal{Y}^0(M)$ & Yes & No & No & Yes\\
\hline
$\mathcal{Y}^-(M)$ & No & Yes & No & Yes\\
\end{tabular}
\end{center}
\vspace{.05in}

\noindent
In the class ($\mathcal{Y}^0(M), \| \sigma\|^2 \equiv 0, \tau =0$), any constant is a solution.  For data in all other classes for which solutions exist, the solution is unique.

Since the solvability of (\ref {eq:lich}) has been completely determined in the constant mean curvature case, the next question is whether the set of solutions has an interesting and useful mathematical structure \cite{isen2}.  We study this question in this paper.  Since any constant is a solution in the class ($\mathcal{Y}^0(M), \| \sigma\|^2 \equiv 0, \tau =0$), we cannot expect to have any compactness of solutions.  In fact, let $c$ be any positive constant, from (\ref {eq:lich}) we know that the function $c\phi$ satisfies $$\Delta (c\phi)  =  \frac{1}{8} R (c\phi) - \frac{1}{8}\|c^4\sigma\|^2(c\phi)^{-7}+\frac{1}{12}(c^{-2}\tau)^2(c\phi)^5.$$When $c \to 0$, $\|c^4\sigma\|^2 \to 0$ and $(c^{-2}\tau)^2 \to \infty$; when $c \to \infty$, $\|c^4\sigma\|^2 \to \infty$ and $(c^{-2}\tau)^2 \to 0.$  So in order to get uniform $C^0$ estimates for $\phi$, we will need to control $\|\sigma\|^2$ and $\tau^2$ from above and below as well.  The following theorem shows that with an additional $C^{0, \alpha}$ bound on $\|\sigma\|^2$, in the remaining five classes we do have compactness of solutions to (\ref {eq:lich}) with respect to $\sigma$ and $\tau$.  For simplicity in the rest of the paper we denote $\|\sigma\|^2$ as $\sigma ^2$.

\begin{thm}
\label {thm:main}
Let $(\lambda _{ij}, \sigma_{ij}, \tau)$ be a determined data on $M^3$, where $\tau$ is a constant.  Let $\phi$ be a positive solution of the Lichnerowicz equation $(\ref {eq:lich}).$

\begin{itemize}
\item In classes $\left(\mathcal{Y}^+, \sigma ^2 \not\equiv0, \tau \neq 0 \right ), \left(\mathcal{Y}^0, \sigma ^2 \not\equiv0, \tau \neq 0 \right ),$ and $ (\mathcal{Y}^-, \sigma ^2 \not\equiv0, \tau \neq 0  ),$ 
assume $C _1^{-1} < \sigma ^2 < C_1 $, $\| \sigma^2\|_{C^{0,\alpha}(M)} <C_1$, and $C_2^{-1} < \tau ^2 < C_2 $ for some constants $C_1, C_2>0$ and $0<\alpha <1.$  Then there exists a constant $C=C(C_1,C_2,\lambda)>0$ such that $C^{-1} \leq  \phi  \leq C$ on $M$.

\item In class $\left(\mathcal{Y}^+, \sigma ^2 \not\equiv0, \tau = 0 \right ),$ assume $C _1^{-1} < \sigma ^2 < C_1 $ and $\| \sigma^2\|_{C^{0,\alpha}(M)} <C_1$ for some constants $C_1>0$ and $0<\alpha <1.$  Then there exists a constant $C=C(C_1,\lambda)>0$ such that $C^{-1} \leq  \phi  \leq C$ on $M$.

\item  In class $\left(\mathcal{Y}^-, \sigma ^2 \equiv0, \tau \neq 0 \right ),$ assume $C_2^{-1} < \tau ^2 < C_2 $ for some constant $C_2>0.$  Then there exists a constant $C=C(C_2,\lambda)>0$ such that $C^{-1} \leq  \phi  \leq C$ on $M$.

\end{itemize}

\end{thm}

\noindent
By a bootstrap argument which will be carried out in detail in the proof of this theorem, these estimates imply that $\phi$ has a uniform $C^3$ norm bound and therefore a sequence of solutions will produce another solution in its limit.  We prove these estimates in the next three sections.


\section  {The Classes $\left(\mathcal{Y}^+, \sigma ^2 \not\equiv0, \tau \neq 0 \right ), \left(\mathcal{Y}^0, \sigma ^2 \not\equiv0, \tau \neq 0 \right ),$ and $ \left(\mathcal{Y}^-, \sigma ^2 \not\equiv0, \tau \neq 0 \right )$}
\label {section:case1}

First,{ \it the lower bound}.

Consider the $1$-variable function $f(t)=\frac{1}{12}\tau^2 t ^3 + \frac{1}{8}Rt^2-\frac{1}{8}\sigma^2$.  Since $f(0)=-\frac{1}{8}\sigma^2 < -\frac{1}{8}C_1^{-1}$ and $f$ is continuous, there exists a constant $\epsilon=\epsilon(C_1, C_2, R(\lambda))>0$ such that if $|t|<\epsilon$ then $f(t)<0$.  Therefore if $0<\phi<\epsilon^{\frac{1}{4}}$, then $\frac{1}{12}\tau ^2\phi ^{12} + \frac{1}{8}R\phi^8-\frac{1}{8}\sigma^2<0$, which implies that $\Delta \phi = \frac{1}{8}R\phi - \frac{1}{8}\sigma^2 \phi ^{-7} + \frac{1}{12}\tau^2 \phi ^5 <0 $.  Since $\Delta \phi$ cannot be negative at the minimum point of $\phi$, we then know that the minimum of $\phi$ cannot be smaller than $\epsilon^{\frac{1}{4}}$.  This proves the lower bound $\phi \geq \epsilon^{\frac{1}{4}}$.

Next, {\it the upper bound.}

Suppose there is no uniform upper bound on $\phi$, then we can find sequences $\{ \phi _i \},$ $ \{ \sigma _i \},$ $ \{ \tau _i \},$ and $\{ x _i \}$, such that   
\begin{equation*}
\Delta \phi _i = \frac{1}{8}R\phi _i - \frac{1}{8}\sigma _i ^2 \phi _i ^{-7} + \frac{1}{12}\tau_i ^2 \phi _i ^5, 
\end{equation*}
\noindent
where $\phi _i >0$, $C _1^{-1} < \sigma _i^2 < C_1 $, $\| \sigma _i^2\|_{C^{0,\alpha}(M)} <C_1$, $C_2^{-1} < \tau_i ^2 < C_2 $, and $\displaystyle \max _{M} \phi _i = \phi _i(x_i) \to \infty$.

Let $x=(x^1, x^2, x^3)$ be the geodesic normal coordinates centered at each of the points $x_i$, and let $y=\phi _i^2(x_i) x$. Define $$u_i(y)=\frac{\phi _i\left ( \frac {y}{\phi^2 _i(x_i)}  \right ) }{\phi _i(x_i)}.$$  Then $u_i$ satisfies 

\begin{eqnarray}
\label{eq:case1-after-scaling}
& & \Delta _{\lambda^{(i)}(y)} u_i(y)\\ \nonumber
& = &\frac{1}{8} R \left ( \frac {y}{\phi^2 _i(x_i)}  \right ) u_i(y) \phi _i^{-4}(x_i) 
- \frac{1}{8}\sigma _i ^2  \left ( \frac {y}{\phi^2 _i(x_i)}  \right )  \phi _i ^{-7}  \left ( \frac {y}{\phi^2 _i(x_i)}  \right ) \phi _i^{-5}(x_i) 
+ \frac{1}{12}\tau_i ^2 u _i ^5(y),
\end{eqnarray}

\noindent
where the metric $\lambda^{(i)}(y)=\sum_{k,l=1}^3
\lambda_{kl} \left ( \frac {y}{\phi^2 _i(x_i)}  \right )dy^kdy^l.$
By the definition we also know that $0<u_i \leq 1$ and $u_i(0)=1$.  On any compact subset $\Omega$ of the $y$-plane, since we have proved that $\phi$ has a positive lower bound, the right hand side of (\ref {eq:case1-after-scaling}) has a uniformly bounded $L^p$ norm for arbitrary $p>1$.  Then by standard elliptic estimates $ u_i$ has a uniform $W^{2,p}$ norm bound, which by the Sobolev embedding theorem implies that $u_i$ is uniformly bounded in the $C^{1,\alpha}$ norm when $p$ is large enough such that $1-\frac{3}{p} \geq \alpha$.  

\noindent
Next we show that $- \frac{1}{8}\sigma _i ^2  \left ( \frac {y}{\phi^2 _i(x_i)}  \right )  \phi _i ^{-7}  \left ( \frac {y}{\phi^2 _i(x_i)}  \right ) \phi _i^{-5}(x_i)$, i.e. the second term on the right hand side of (\ref {eq:case1-after-scaling}), has bounded $C^{0,\alpha}$ norm on $\Omega.$  Given the assumption on $\sigma _i^2$ and the fact that $\phi_i$ is bounded below and $\phi _i(x_i) \to \infty$, we only need to obtain a $C^{0,\alpha}$ bound on $\phi _i ^{-7}  \left ( \frac {y}{\phi^2 _i(x_i)}  \right ) \phi _i^{-5}(x_i)$.  Let $y_1$ and $y_2$ be any two points on $\Omega,$ by the definition of $u_i$  
\begin{eqnarray}
\label{eq:holder}
\frac{\Big |\phi _i ^{-7}  \left ( \frac {y_1}{\phi^2 _i(x_i)}  \right ) -\phi _i ^{-7}  \left ( \frac {y_2}{\phi^2 _i(x_i)}  \right ) \Big | }{\phi _i^5(x_i)|y_1-y_2|^{\alpha}} & = & \frac{\big | \phi _i^{-7}(x_i)\big (u_i^{-7}(y_1)-u_i^{-7}(y_2)\big ) \big |}{\phi _i^5(x_i)|y_1-y_2|^{\alpha}} \nonumber  \\
& = & \frac{\big | u_i^{-7}(y_1)-u_i^{-7}(y_2)\big |}{\phi _i^{12}(x_i)|y_1-y_2|^{\alpha}}.
\end{eqnarray}

\noindent
Note that 
\begin{eqnarray*}
\big | u_i^{-7}(y_1)-u_i^{-7}(y_2)\big | & = &  \Bigg | \int _{u_i(y_2)}^{u_i(y_1)} \frac{d}{dt}(t^{-7}) dt \Bigg | \\
& = & \Bigg | 7 \int _{u_i(y_2)}^{u_i(y_1)} t^{-8} dt \Bigg | \\
& \leq & 7 \left (\max _{y \in \Omega} u_i^{-8}(y) \right ) |u_i(y_1)-u_i(y_2) | \\
& = & 7 \left (\max _{y \in \Omega} \phi _i^{-8} \left ( \frac {y}{\phi^2 _i(x_i)}  \right ) \phi _i^8(x_i)  \right ) |u_i(y_1)-u_i(y_2) | \\
& \leq & C \phi _i^8(x_i)|u_i(y_1)-u_i(y_2) | 
\end{eqnarray*}

\noindent 
for some constant $C$, where the last inequality has used the positive lower bound on $\phi$.  Therefore 
\begin{eqnarray*}
 \frac{\big | u_i^{-7}(y_1)-u_i^{-7}(y_2)\big |}{\phi _i^{12}(x_i)|y_1-y_2|^{\alpha}} & \leq & \frac {C \phi _i^8(x_i)|u_i(y_1)-u_i(y_2) | }{\phi _i^{12}(x_i)|y_1-y_2|^{\alpha}}\\
& \leq & \frac{C\|u_i\|_{C^{0, \alpha}(\Omega)}}{\phi_i^4(x_i)}\\
& \leq &  C(C_1, C_2, \Omega).
\end{eqnarray*}

\noindent
By (\ref {eq:holder}) this gives a uniform bound on the $C^{0,\alpha}$ norm of $\phi _i ^{-7}  \left ( \frac {y}{\phi^2 _i(x_i)}  \right ) \phi _i^{-5}(x_i)$ on $\Omega$.  As we explained earlier, this then implies a uniform $C^{0,\alpha}$ bound on the second term of the right hand side of (\ref {eq:case1-after-scaling}).  Then since the other two terms on that side both have $C^{1,\alpha}$ bound, by the Schauder estimates we have uniform $C^{2,\alpha}$ bound on $u_i$ on the compact set $\Omega$.  This implies that a sequence of $u_i$ converges in $C^2$ norm to some function $u$ on $\Omega$ where $u$ satisfies $0\leq u \leq 1$ and $u(0)=1$.  By the assumptions on $\sigma _i^2$ and $\tau _i^2$ we also know that on $\Omega$, passing to subsequences $\{ \sigma _i ^2\}$ converges to some function $C_1^{-1} \leq \sigma ^2 \leq C_1$ and $\{\tau _i ^2 \}$ converges to some constant $C_2^{-1} \leq \tau ^2 \leq C_2$.  Additionally, the metrics $\lambda^{(i)}$ converge to the Euclidean metric on $\Omega$.  Then we let $i \to \infty$ on both sides of (\ref {eq:case1-after-scaling}) for $y \in \Omega$. Since $\phi_i(x_i) \to \infty$ and $R$ and $ \phi_i^{-7}$ are bounded above, in the limit the equation becomes 
\begin{equation}
\label {eq:liou}
\Delta _{\delta} u = \frac{1}{12} \tau^2 u^5
\end{equation}

\noindent
where $\Delta _{\delta}$ denotes the Euclidean Laplacian.  Since $\Omega$ is arbitrary, we thus have obtained a function $0\leq u \leq 1, u(0)=1$ which satisfies (\ref {eq:liou}) on $\mathbf{R}^3(y)$.  This is in fact impossible due to the following lemma, and therefore we have completed the proof in these three classes.

\begin{lemma}
\label {lemma:liou}
There does not exist any function $0\leq u \leq 1, u(0)=1$ which satisfies equation $(\ref {eq:liou})$ on $\mathbf{R}^3$.
\end{lemma}

\pf To simplify the proof we assume the constant $\frac{1}{12} \tau ^2=1$.  We will use the moving plane method as in \cite{Lin} by C.S. Lin, but in our case the proof is much easier with the extra assumptions on $u$.

Suppose such a function $u$ exists.  Then for any $t \in \mathbf{R}$ and any $y=(y^1,y^2, y^3)$, denote $y_t=(2t-y^1, y^2, y^3)$ and define $u_t(y)=u(y_t)$.  Then $u_t$ also satisfies $\Delta _{\delta} u_t=u_t^5$.  We claim that for any $t$, $u_t(y) \geq u(y)$ whenever $y^1<t$.  

\noindent
Suppose this is not true.  Define the function $w_t(y)=u_t(y)-u(y)$, then for some $t$ we have $$\displaystyle \inf _{\{y:y^1< t\}}w_t(y)<0.$$  Now let $$g_t(y)= \ln \left ( (t-y^1+2)^2+(y^2)^2\right)+\ln \left ( (t-y^1+2)^2+(y^3)^2\right)$$ for $y\in \mathbf{R}^3 $ with $y^1 < t$.  By this definition we know that 
$$g_t(y) > \ln 2^2 >1, \hspace {.2in} \lim _{|y| \to \infty} g_t(y) = \infty, \hspace{.2in} \text {and} \hspace{.2in} \Delta _{\delta} g_t=0.$$  
Define $$\bar{w}_t(y)=\frac{w_t(y)}{g_t(y)},$$ then  $\displaystyle \inf _{\{y:y^1< t\}}\bar{w}_t(y)<0$ because $\displaystyle \inf _{\{y:y^1< t\}}w_t(y)<0$ and $g_t(y)>1$.  Because $\displaystyle\lim _{|y| \to \infty} g_t(y) = \infty$ and $w_t$ is bounded, we have $\displaystyle\lim _{|y| \to \infty} \bar{w}_t(y) = 0.$  Therefore $ \displaystyle \inf_{\{y:y^1< t\}}\bar{w}_t(y)$ is achieved at some point $y_0$, and $\bar{w}_t(y_0)<0.$

\noindent
On one hand, from $g_t(y_0)>1$ and $\bar{w}_t(y_0)<0$ we know that $w_t(y_0)<0$, i.e. $u_t(y_0) < u(y_0)$.  This implies that
\begin{eqnarray*}
\Delta_{\delta} w_t(y_0) & = & \Delta_{\delta} u_t(y_0)-\Delta_{\delta} u(y_0)\\
& = & u_t^5(y_0)-u^5(y_0)\\
& < & 0.
\end{eqnarray*}

\noindent
On the other hand, since $y_0$ is a minimum point for $\bar{w}_t$, we have that $\nabla _{\delta}\bar{w}_t(y_0)=0$ and $\Delta_{\delta} \bar{w}_t(y_0)\geq 0$.  Combined with $g_t>1$ and $\Delta_{\delta} g_t=0 $, this leads to
\begin{eqnarray*}
\Delta_{\delta} w_t(y_0) & = & \Delta_{\delta} (\bar{w}_t g_t  )  (y_0)\\
& = & g_t(y_0) \Delta_{\delta} \bar{w}_t(y_0) + \bar{w}_t(y_0) \Delta_{\delta} g_t(y_0) + 2 \nabla _{\delta}\bar{w}_t(y_0) \cdot \nabla _{\delta}g_t(y_0)   \\
& = &  g_t(y_0) \Delta_{\delta} \bar{w}_t(y_0)\\
& \geq & 0.
\end{eqnarray*}

\noindent
Thus we have reached a contradiction.  Therefore $u_t(y)\geq u(y)$ for any $t\in \mathbf{R}$ and $y^1<t$.  

Then because $u(0)=1$ is the maximum, $u$ must be identically equally to $1$ on the positive $y^1$-axis.  Since the equation $\Delta _{\delta} u = u^5$ is invariant under rotations about the origin, this implies that $u$ is identically equal to $1$ on any half line emanating from the origin.  Therefore we know that $u \equiv 1$.  However, this contradicts the equation  $\Delta _{\delta} u = u^5,$ and the proof of this lemma is finished. 

\stop

\section{ The Class $\left(\mathcal{Y}^+, \sigma ^2 \not\equiv0, \tau = 0 \right )$}
\label {section:case 2}

In this class the Lichnerowicz equation becomes 
$$\Delta \phi = \frac{1}{8}R\phi - \frac{1}{8} \sigma ^2 \phi ^{-7}.$$  Due to the conformal invariant property of this equation, we can assume $R$ to be a positive constant.

First {\it the lower bound.}

If $\phi < \left ( \frac{\sigma^2}{R} \right )^{\frac{1}{8}}$, then $\Delta \phi = \frac{1}{8}R\phi - \frac{1}{8} \sigma ^2 \phi ^{-7}<0$, hence $\phi$ cannot reach a minimum.  Therefore  $\phi \geq  \left ( \frac{\sigma^2}{R} \right )^{\frac{1}{8}} \geq \left ( C_1^{-1}R^{-1} \right )^{\frac{1}{8}}.$ 

Next {\it the upper bound.}

Suppose there is no uniform upper bound on $\phi$, then we can find sequences $\{ \phi _i \}, \{ \sigma _i \},$ and $\{ x _i \}$, such that   
\begin{equation*}
\Delta \phi _i = \frac{1}{8}R\phi _i - \frac{1}{8}\sigma _i ^2 \phi _i ^{-7}, 
\end{equation*}
\noindent
where $\phi _i >0$, $C _1^{-1} < \sigma _i^2 < C_1 $, $\| \sigma _i^2\|_{C^{0,\alpha}(M)} <C_1$, and $\displaystyle \max _{M} \phi _i = \phi _i(x_i) \to \infty$.

\noindent
Now define another sequence of functions $\{v_i\}$ on $M$ by $$v_i(x)=\frac{\phi_i(x)}{\phi _i(x_i)}.$$  Then $0<v_i(x) \leq 1$, $v_i(x_i) = 1$, and 
\begin{equation}
\label {eq:case2-after-scaling}
\Delta v_i = \frac{1}{8}R v_i - \frac{1}{8}\sigma_i^2 \phi _i^{-7}\phi _i^{-1}(x_i).
\end{equation}

Since we already proved the positive lower bound on $\phi_i$, the right hand side of (\ref {eq:case2-after-scaling}) has uniformly bounded $L^p(M)$ norm for any $p>1$.  Then by standard elliptic estimates $v_i$ has bounded $W^{2,p}$ norm, which leads to uniform $C^{1, \alpha}$ bound when $p$ is large enough.  Next we show that $\phi _i^{-7}\phi _i^{-1}(x_i)$ has uniform $C^{0,\alpha}$ bound as follows. 

\noindent
For any two points $p_1, p_2 \in M$, denote $d_{\lambda}(p_1,p_2)$ as $|p_1-p_2|$.  Then 

\begin{eqnarray*}
 \frac{\big | \phi_i^{-7}(p_1)-\phi_i^{-7}(p_2)\big |}{\phi _i(x_i)|p_1-p_2|^{\alpha}} & = &  \frac{\big | v_i^{-7}(p_1)-v_i^{-7}(p_2)\big |}{\phi ^8 _i(x_i)|p_1-p_2|^{\alpha}}\\
& = &  \frac{\Big | \int _{v_i(p_2)}^{v_i(p_1)} \frac{d}{dt}(t^{-7}) dt \Big |}{\phi ^8 _i(x_i)|p_1-p_2|^{\alpha} }  \\
& \leq & \frac{ 7 \left (\max _{M} v_i^{-8} \right ) |v_i(p_1)-v_i(p_2) |} {\phi ^8 _i(x_i)|p_1-p_2|^{\alpha}  }  \\
& = & \frac{7 \left (\max _{M} \phi _i^{-8}   \phi _i^8(x_i)  \right ) |v_i(p_1)-v_i(p_2) | }{\phi ^8 _i(x_i)|p_1-p_2|^{\alpha}  }\\
& \leq & 7 \left ( \max _{M}  \phi _i^{-8} \right ) \|v_i \| _{C^{0,\alpha}(M)}\\
& \leq & C.
\end{eqnarray*}

\noindent
Then since $\sigma _i^2$ is also bounded in $C^{0,\alpha}$ norm and $v_i$ has bounded $C^{1, \alpha}$ norm, the right hand side of (\ref {eq:case2-after-scaling}) has bounded $C^{0,\alpha}$ norm.  This implies a uniform $C^{2,\alpha}$ bound on $v_i$ by the Schauder estimates.  

Thus a subsequence of $\{v_i\}$ converges in the $C^2$ norm to some function $v$.  Since $\sigma_i^2$ and $\phi _i^{-7}$ are both bounded above, we can take the limits of both sides of (\ref {eq:case2-after-scaling}) as $i \to \infty$ and obtain $\Delta v = \frac{1}{8}Rv$ on $M,$ where $0 \leq v \leq 1$.  Because $R>0$, the maximal principle implies $v \equiv 0$.  However, since $M$ is compact, there exists a point $x_0 \in M$ such that passing to a subsequence $x_0=\displaystyle \lim_{i \to \infty} x_i.$  Then $v(x_0) = \displaystyle \lim_{i \to \infty} v_i(x_i) =1.$  This contradiction finishes the proof of the upper bound.

\section{ The Class $\left(\mathcal{Y}^-, \sigma ^2 \equiv0, \tau \neq 0 \right )$}
\label {section:case 3}

In this class the Lichnerowicz equation becomes 
$$\Delta \phi = \frac{1}{8}R\phi +\frac{1}{12} \tau ^2 \phi ^5.$$  Due to the conformal invariant property of this equation, we can assume $R$ to be a negative constant.

First {\it the lower bound.}

If $\phi < \left ( \frac{-3R}{2\tau^2} \right )^{\frac{1}{4}}$, then $\Delta \phi = \frac{1}{8}R\phi + \frac{1}{12} \tau ^2 \phi ^{5}<0$, hence $\phi$ cannot reach a minimum.  Therefore  $\phi \geq \left ( \frac{-3R}{2\tau^2} \right )^{\frac{1}{4}}   \geq ( -3R )^{\frac{1}{4}}( 2^{-1}C_2^{-1})^{\frac{1}{4}}.$ 

Next {\it the upper bound.}

Suppose there is no uniform upper bound on $\phi$, then we can find sequences $\{ \phi _i \},$ $ \{ \tau _i \},$ and $\{ x _i \}$, such that   
\begin{equation*}
\Delta \phi _i = \frac{1}{8}R\phi _i  + \frac{1}{12}\tau_i ^2 \phi _i ^5, 
\end{equation*}
\noindent
where $\phi _i >0$, $C_2^{-1} < \tau_i ^2 < C_2 $, and $\displaystyle \max _{M} \phi _i = \phi _i(x_i) \to \infty$.

Let $x=(x^1, x^2, x^3)$ be the geodesic normal coordinates centered at each of the points $x_i$, and let $y=\phi _i^2(x_i) x$. Define $u_i(y)=\frac{\phi _i\left ( \frac {y}{\phi^2 _i(x_i)}  \right ) }{\phi _i(x_i)}.$   Then by the same argument as in Section \ref{section:case1}, we can show that on any compact subset of $\mathbf{R}^3(y)$, a subsequence of $\{u_i\}$ converges in $C^2$ norm to a function $u$ which satisfies $0 \leq u \leq 1$, $u(0)=1$ and $$\Delta _{\delta} u = \frac{1}{12}\tau^2 u^5$$ on $\mathbf{R}^3,$ where $\Delta _{\delta}$ is the Euclidean Laplacian and $\tau^2 = \displaystyle \lim _{i \to \infty} \tau_i^2$.  In fact, the proof of this convergence is easier in this class than in Section  \ref{section:case1}, because here we do not have the term $\frac{1}{8}\sigma _i^2\phi _i^{-7}$ in the Lichnerowicz equations.  Finally, by Lemma \ref {lemma:liou} such a function $u$ does not exist, and this contradiction completes the proof.

\bibliographystyle{plain}
    \bibliography{thesis}

\end{document}